\documentclass[graybox]{svmult}

\usepackage{type1cm}
\usepackage{makeidx}
\usepackage{multicol}
\usepackage[bottom]{footmisc}
\usepackage{newtxtext}
\usepackage[varvw]{newtxmath}

\makeindex

\usepackage[dvipsnames]{xcolor}

\usepackage[utf8]{inputenc}
\usepackage[T1]{fontenc}

\usepackage{stmaryrd}

\usepackage[style=numeric,backend=biber]{biblatex}

\usepackage{tikz, tikz-cd}
\usetikzlibrary{calc}

\usepackage{quiver}

\usepackage{graphicx}

\usepackage{proof}

\usepackage{hyperref}

\usepackage{todonotes}

\usepackage{float}

\newtheorem*{thm*}{Theorem}

\theoremstyle{definition}
\newtheorem{defn/}{Definition}
\newtheorem*{defn*/}{Definition}

\newcommand{\defnendsymbol}%
{%
  \mathbin{\rotatebox[origin=c]{-45}{$\parallel$}}%
}

\newenvironment{defn}
  {\renewcommand{\qedsymbol}{$\defnendsymbol$}%
   \pushQED{\qed}\begin{defn/}}
  {\popQED\end{defn/}}

\newenvironment{defn*}
  {\renewcommand{\qedsymbol}{$\defnendsymbol$}%
   \pushQED{\qed}\begin{defn*/}}
  {\popQED\end{defn*/}}

\newcommand*{\myimportant}[1]{\textcolor{MidnightBlue}{#1}}

\IfFileExists{../preamble.tex}{\input{../preamble.tex}}{}

\usepackage[scr=boondox,  
            cal=esstix]   
           {mathalpha}

\usepackage{float}
\usepackage{nicematrix}
\usepackage{multirow}

\usetikzlibrary{decorations.markings, arrows.meta, calc, intersections,bending}

\newcommand{\Fuk}{\mathcal{Fuk}}

\author{Christina Grossack}
\title{Explicitly Computing with Fukaya Categories of Surfaces with Boundary}

\begin{document}
\maketitle

\abstract{Fukaya categories are deep and rich invariants of symplectic 
manifolds which are notoriously difficult to compute explicitly. In the 
case of surfaces, however, the situation is simple, combinatorial,
and is very well understood (at least by experts). In this expository paper 
we will give an introduction with many examples to welcome newcomers to the 
area and hopefully equip them with the tools to independently compute Fukaya 
categories of surfaces.}

\section{Introduction}

Fukaya categories are rich invariants of symplectic manifolds which count 
intersection points of Lagrangian submanifolds, and
are interesting from many perspectives. Historical 
motivation comes from physics and mirror symmetry~\cite{kontsevichHomologicalAlgebraMirror1994,
bocklandtGentleIntroductionHomological2021} and more recently the case of 
surfaces was directly related to the 
well studied (locally) gentle algebras from representation 
theory~\cite{haidenFlatSurfacesStability2017, opperGeometricModelDerived2018,
lekiliDerivedEquivalencesGentle2020}. To a general symplectic manifold $M$, we 
want to assign a category\footnote{Really an $A_\infty$-category, where 
composition need only be associative up to coherent higher homotopy} 
whose objects are Lagrangian submanifolds (possibly equipped with bonus data 
such as a local system, a spin structure, a Morse function, etc.) and arrows 
between two submanifolds $L_1$ and $L_2$ are points in the (transverse) intersection 
$L_1 \cap L_2$. Composition (and higher $A_\infty$-operations) are computed 
using the space of solutions to a differential equation whose boundary 
conditions are given by the Lagrangians and the intersection points one wants
to compose. This is obviously difficult to compute in practice, and there are 
many variants of the Fukaya category which have different computational 
properties. See, for instance,~\cite[Chapter 6]{bocklandtGentleIntroductionHomological2021}
and~\cite[Section 2.3]{aurouxBeginnersIntroductionFukaya2014}.
Luckily, in the case of surfaces with nonempty boundary, the situation 
simplifies considerably. 
In this case Lagrangians are ``just'' curves in the surface, and the choice of 
symplectic structure does not play a role\footnote{note though, that the
symplectic structure \emph{is} relevant for closed surfaces}. Moreover, in 
many situations one can arrange for the Fukaya category to be a dg-category, 
where the higher the higher $A_\infty$-operations vanish. This may feel more 
concrete to many readers\footnote{it certainly felt more concrete to the author 
when she was first entering the subject}. 

Thus the computation of the 
\myimportant{Partially Wrapped Fukaya Category} (sometimes called the 
\myimportant{Topological Fukaya Category} in this context, to emphasize that 
it does not depend on a choice of symplectic form) becomes pleasantly 
combinatorial, though the author is not currently aware of any places in the 
literature where multiple ``large'' concrete examples are worked out in detail.
The author recently gave a talk showing how to compute 
examples in practice, and was encouraged to record these examples, which 
led to the current expository article.

Nothing in this article is original except for the exposition, 
which we hope can serve as a shorter introduction to the computational details 
of Fukaya categories of surfaces, aimed at newcomers to the field.
For important and interesting theoretical details about \emph{why} these 
computations work, we again point the reader to the foundational 
papers~\cite{haidenFlatSurfacesStability2017, opperGeometricModelDerived2018, 
lekiliDerivedEquivalencesGentle2020} and the recent 
textbook~\cite{bocklandtGentleIntroductionHomological2021} which 
provides lots of important context and motivation related to mirror symmetry,
in addition to other techniques which we had to cut from this note for 
reasons of space\footnote{Such as information about $A_\infty$-structures and 
gluing Fukaya categories by taking (homotopy) colimits}.
For the reader who enjoys seeing examples worked out in real time, we highly 
recommend the lecture series of 
Claire Amiot~\cite{institutfourierClaireAmiotFukaya2024b,
institutfourierClaireAmiotFukaya2024c,
institutfourierClaireAmiotFukaya2024a,
institutfourierClaireAmiotFukaya2024} and 
Sybille Schroll~\cite{hausdorffcenterformathematicsGeometricModelBounded2020b,
hausdorffcenterformathematicsGeometricModelBounded2020a,
hausdorffcenterformathematicsGeometricModelBounded2020}, both of which are 
available on Youtube.

\begin{acknowledgement}
The author would like to thank the organizers of the 
Women and Gender Minorities in Symplectic and Contact Geometry and Topology
(WiSCon) special session -- Bahar Acu, Catherine Cannizzo, Sierra Knavel, 
and Morgan Weiler -- both for the opportunity to speak and for the incredible 
community they're fostering. Additionally, she would like to thank 
Peter Samuelson and Shane Rankin for many helpful conversations while 
learning this material, and especially Catherine Cannizzo \emph{again} for 
answering many questions and helping the author organize a learning seminar 
on this topic. Finally, she would like to thank the AWM and the NSF for their 
financial support.
\end{acknowledgement}

\section{Marked Surfaces and Dissections}

\begin{defn}\label{def:MarkedSurface}
  A \myimportant{Marked Surface} is an oriented topological surface $S$ 
  with nonempty boundary $\partial S$ and a finite set of 
  marked points $M \subseteq \partial S$ called \myimportant{Stops}.
\end{defn}

We think of our boundary components as extending to infinity away from the 
stops, as shown in Figure~\ref{fig:MarkedSurfaces}.

\begin{figure}
    \centering
    
    \begin{minipage}{0.15\textwidth}
        \centering
        \resizebox{\linewidth}{!}%
        {%
          \begin{tikzpicture}[line cap=round, line join=round]
            
            \def\outerR{2.0} 
            \def\innerR{0.45} 
            \def\dotR{0.12}  
            
            \def\angleL{225} 
            \def\angleR{315} 
            
            \draw[thick] (0,0) circle [radius=\outerR];
            \draw[thick] (0,0) circle [radius=\innerR];

            \fill (90:\outerR) circle [radius=\dotR];  
            \fill (\angleL:\outerR) circle [radius=\dotR]; 
            \fill (\angleR:\outerR) circle [radius=\dotR]; 
          \end{tikzpicture}
        }
    \end{minipage}%
    \quad
    \begin{minipage}{0.15\textwidth}
        \centering
        \resizebox{\linewidth}{!}%
        {%
          \begin{tikzpicture}[line cap=round, line join=round]
              
              \def\R{2.0}      
              \def\ellE{0.4}   
              \def\h{4.5}      
              
              \def\slotW{0.48} 
              \def\slotD{2.2}  
              \def\topY{2.0}   
              
              \pgfmathsetmacro{\xLeft}{-\R*cos(45)}  
              \pgfmathsetmacro{\xRight}{\R*cos(45)}  

              \def\alignTop{3}
              \draw[dashed, thin, black!60] ({-\R}, {\topY}) -- ({-\R}, {\alignTop}); 
              \draw[dashed, thin, black!60] ({\R}, {\topY}) -- ({\R}, {\alignTop});   
              
              \draw[thick] ({-\R}, {-\h/2}) -- ({-\R}, {\topY}); 
              \draw[thick] ({\R}, {-\h/2}) -- ({\R}, {\topY});   
              
              \def\backSlotAngle{83.1}
              \def\backSlotAngleL{96.9}
              
              \draw[thick, black!90] ({\R}, {\topY}) arc [start angle=0, end angle=\backSlotAngle, x radius=\R, y radius=\ellE];
              \draw[thick, black!90] ({-\R*cos(\backSlotAngle)}, {\topY + \ellE*sin(\backSlotAngleL)}) arc [start angle=\backSlotAngleL, end angle=180, x radius=\R, y radius=\ellE];
              
              \draw[thick, dashed, black!70] 
                  ({-\slotW/2}, {\topY + \ellE*sin(\backSlotAngleL)}) 
              -- ({-\slotW/2}, {\topY + \ellE*sin(\backSlotAngleL) - \slotD})
              arc [start angle=180, end angle=360, radius={\slotW/2}]
              -- ({\slotW/2}, {\topY + \ellE*sin(\backSlotAngle)});
              
              \pgfmathsetmacro{\aLout}{360 - acos((\xLeft - \slotW/2)/\R)}
              \pgfmathsetmacro{\aLin}{360 - acos((\xLeft + \slotW/2)/\R)}
              \pgfmathsetmacro{\aRin}{360 - acos((\xRight - \slotW/2)/\R)}
              \pgfmathsetmacro{\aRout}{360 - acos((\xRight + \slotW/2)/\R)}

              \draw[thick] ({-\R}, {\topY}) arc [start angle=180, end angle=\aLout, x radius=\R, y radius=\ellE];
              \draw[thick] ({\R*cos(\aLin)}, {\topY + \ellE*sin(\aLin)}) arc [start angle=\aLin, end angle=\aRin, x radius=\R, y radius=\ellE];
              \draw[thick] ({\R*cos(\aRout)}, {\topY + \ellE*sin(\aRout)}) arc [start angle=\aRout, end angle=360, x radius=\R, y radius=\ellE];

              \draw[thick] ({\xLeft - \slotW/2}, {\topY + \ellE*sin(\aLout)}) 
                        -- ({\xLeft - \slotW/2}, {\topY - \slotD})
                      arc [start angle=180, end angle=360, radius={\slotW/2}]
                        -- ({\xLeft + \slotW/2}, {\topY + \ellE*sin(\aLin)});

              \draw[thick] ({\xRight - \slotW/2}, {\topY + \ellE*sin(\aRin)}) 
                        -- ({\xRight - \slotW/2}, {\topY - \slotD})
                      arc [start angle=180, end angle=360, radius={\slotW/2}]
                        -- ({\xRight + \slotW/2}, {\topY + \ellE*sin(\aRout)});

              \draw[thick] (0, {-\h/2}) ellipse [x radius=\R, y radius=\ellE];

              \def\alignBot{-3}
              \def\botStart{-2}
              \draw[dashed, thin, black!60] ({-\R}, {\botStart}) -- ({-\R}, {\alignBot}); 
              \draw[dashed, thin, black!60] ({\R}, {\botStart}) -- ({\R}, {\alignBot});   
            \end{tikzpicture} 
        }
    \end{minipage}
    \hfill 
    \begin{minipage}{0.15\textwidth}
        \centering
        \resizebox{\linewidth}{!}%
        {%
          \begin{tikzpicture}
          
          \draw (0, 2) ellipse (1.5 and 0.4);

          \draw (-1.5, 2) .. controls (-1.5, -0.5) and (-0.9, -1.5) .. (0, -1.5)
                          .. controls (0.9, -1.5) and (1.5, -0.5) .. (1.5, 2);

          \fill[black!80] (1.3, 1.8) circle (0.10);

          \draw (0.1, 1.1) 
            .. controls (-0.1, 0.8) and (-0.2, 0.6) .. (0.0, 0.5)    
            .. controls (0.2, 0.4) and (0.2, 0.8) .. (-0.05, 0.75)   
            .. controls (-0.15, 0.7) and (0.0, 0.4) .. (0.2, 0.2);   

          \end{tikzpicture}
        }
    \end{minipage}
    \quad
    \begin{minipage}{0.15\textwidth}
        \centering
        \resizebox{\linewidth}{!}%
        {%
          \begin{tikzpicture}[
              thick,
              >={Stealth[length=6pt,width=6pt]},
              single arrow/.style={%
                  decoration={markings, mark=at position 0.55 with {\arrow{>}%
                  }},
                  postaction={decorate}
              },
              double arrow/.style={
                  decoration={markings, 
                      mark=at position 0.51 with {\arrow{>}},
                      mark=at position 0.59 with {\arrow{>}}
                  },
                  postaction={decorate}
              }
          ]

          \coordinate (A) at (0,0);
          \coordinate (B) at (4,0);
          \coordinate (C) at (4,4);
          \coordinate (D) at (0,4);

          \draw[double arrow] (A) -- (B); 
          \draw[single arrow] (B) -- (C); 
          \draw[double arrow] (D) -- (C); 
          \draw[single arrow] (A) -- (D); 

          \draw (0,0) .. controls (0.2, 1.5) and (1.5, 0.2) .. (0,0) node[pos=0.5, circle, fill=black, inner sep=2.5pt] {};

          \end{tikzpicture}
        }
    \end{minipage}
    \quad
    \begin{minipage}{0.15\textwidth}
      \centering
      \resizebox{\linewidth}{!}
      {%
        \begin{tikzpicture}[scale=2, line cap=round, line join=round]

            \tikzset{main edge/.style={thick, black}}
            \tikzset{hole edge/.style={thick, black}}

            \draw[thick, gray, dashed] (-1, 1.8) -- (-1, 2.3);
            \draw[thick, gray, dashed] (1, 1.8) -- (1, 2.3);

            \draw[main edge] (-1, 1.8) .. controls (-1, 2.2) and (1, 2.2) .. (1, 1.8);

            \draw[main edge] (-1, 1.8) 
                .. controls (-1, 0.4) and (-0.7, -0.8) .. (0, -0.8)
                .. controls (0.7, -0.8) and (1, 0.4) .. (1, 1.8);

            \draw[main edge, fill=white] (-1, 1.8) 
                .. controls (-0.4, 1.5) and (0.0, 1.6) .. (0.35, 1.6) 
                -- (0.35, 1.1) .. controls (0.35, 0.85) and (0.65, 0.85) .. (0.65, 1.1) 
                -- (0.65, 1.6) .. controls (0.7, 1.6) and (0.9, 1.7) .. (1, 1.8); 

            \draw[main edge] (-0.08, 0.35) .. controls (-0.2, 0.1) and (-0.2, -0.2) .. (-0.08, -0.45);
            \draw[hole edge] (-0.12, 0.15) .. controls (-0.02, 0.05) and (-0.02, -0.15) .. (-0.12, -0.25);

        \end{tikzpicture}
      }
    \end{minipage}
    
    \caption{Marked points (called \myimportant{stops}) and their associated 
    geometric description. 
    Note that the boundary extends cylindrically to infinity except in a 
    small neighborhood of the stops. As usual a genus $g$ surface with 
    boundary can be studied as a $2g$-gon with edges identified,
    but it will be useful later to have the edges we cut along end on a 
    boundary component, as shown here.
    }%
    \label{fig:MarkedSurfaces}%
\end{figure}

In the Fukaya category of a marked surface $S$, objects are curves in 
$S \setminus M$ equipped with orientation data (to be described shortly),
local systems of vector spaces. For 
curves $\alpha$ and $\beta$ the homspace $\text{Hom}^\bullet(\alpha, \beta)$ 
is a chain complex containing, for each intersection point 
$p \in \alpha \pitchfork \beta$, the space of linear maps between the stalks 
$\alpha_p$ and $\beta_p$ of the local systems at the point $p$. These linear 
maps are placed in a homological degree $n_p$ related to the orientation data 
of $\alpha$ and $\beta$ which we will describe shortly.

In case the intersection is not transverse, or in case the curves go off to 
infinity, we perturb $\alpha$ by a Hamiltonian isotopy. In the interior of 
the surface this means we want the signed area between $\alpha$ and its 
perturbation $\alpha^H$ to be $0$ and on the boundary this means we flow 
in the direction of the induced orientation on the boundary 
\emph{until we hit a stop}\footnote{In higher dimensional symplectic manifolds
we flow $\alpha$ along the Reeb vector field on the boundary.}. 
This process on the boundary is called \emph{wrapping} and explains why the 
marked points are called stops. 

\begin{figure}
  \centering

  \begin{minipage}{0.25\textwidth}
    \centering
    \resizebox{\linewidth}{!}%
    {%
      \begin{tikzpicture}[line cap=round, line join=round]

        \def\rx{2.0}   
        \def\ry{0.6}   
        \def\h{4.5}    

        \draw[thick, dashed, gray!80] (\rx, 0) arc [start angle=0, end angle=180, x radius=\rx, y radius=\ry];
        \draw[thick] (-\rx, 0) arc [start angle=180, end angle=360, x radius=\rx, y radius=\ry];

        \draw[thick] (-\rx, 0) -- (-\rx, \h);
        \draw[thick] (\rx, 0) -- (\rx, \h);

        \draw[very thick, orange] 
          ({ \rx*cos(235) }, { \ry*sin(235) }) -- ({ \rx*cos(235) }, { \h + \ry*sin(235) });

        \draw[very thick, cyan!90!black] 
          ({ \rx*cos(295) }, { \ry*sin(295) }) -- ({ \rx*cos(295) }, { \h + \ry*sin(295) });

        \draw[thick] (0, \h) ellipse [x radius=\rx, y radius=\ry];

        \fill[black] ({ \rx*cos(335) }, { \h + \ry*sin(335) }) circle (3.5pt);

      \end{tikzpicture}
    }
  \end{minipage}%
  \quad
  \begin{minipage}{0.25\textwidth}
    \centering
    \resizebox{\linewidth}{!}%
    {%
      \begin{tikzpicture}[line cap=round, line join=round]

      \def\rx{2.0}   
      \def\ry{0.6}   
      \def\h{5.5}    

      \draw[thick, dashed, gray!70] (\rx, 0) arc [start angle=0, end angle=180, x radius=\rx, y radius=\ry];
      \draw[thick, gray!70] (-\rx, 0) arc [start angle=180, end angle=360, x radius=\rx, y radius=\ry];

      \draw[thick, gray!70] (-\rx, 0) -- (-\rx, \h);
      \draw[thick, gray!70] (\rx, 0) -- (\rx, \h);

      \def\angleBlue{280}
      \draw[very thick, cyan!90!black] 
        ({ \rx*cos(\angleBlue) }, { \ry*sin(\angleBlue) }) -- ({ \rx*cos(\angleBlue) }, { \h + \ry*sin(\angleBlue) });

      \draw[thick, gray!70] (0, \h) ellipse [x radius=\rx, y radius=\ry];

      \def\anglePoint{335}
      \fill[black] ({ \rx*cos(\anglePoint) }, { \h + \ry*sin(\anglePoint) }) circle (4pt);

      \def\Nturns{3}              
      \def\tStart{180}            
      \def\tEndHelix{295}         
      \def\tVertical{320}         
      \def\hHelix{3.6}            
      \def\hVertical{4.3}         
      
      \def\dtTotal{(\Nturns*360 + \tEndHelix - \tStart)}
      \def\zHeightFunc#1{\hHelix * (#1 - \tStart) / \dtTotal}
      
      \foreach \k in {0, 1, 2} {
        \def\tTurnA{\tStart + \k*360}
        \draw[very thick, orange, opacity=0.9] plot[domain=\tTurnA:{\tTurnA+180}, samples=20] 
          ({ \rx*cos(\x) }, { \ry*sin(\x) + \zHeightFunc{\x} });
        \draw[very thick, orange, dashed, opacity=0.9] plot[domain={\tTurnA+180}:{\tTurnA+360}, samples=20] 
          ({ \rx*cos(\x) }, { \ry*sin(\x) + \zHeightFunc{\x} });
      }
      
      \def\tEndFullWraps{\tStart + \Nturns*360}
      \draw[very thick, orange, opacity=0.9] plot[domain=\tEndFullWraps:{\tEndFullWraps + (\tEndHelix - \tStart)}, samples=20] 
        ({ \rx*cos(\x) }, { \ry*sin(\x) + \zHeightFunc{\x} });

      \coordinate (HelixEnd) at ({ \rx*cos(\tEndHelix) }, { \ry*sin(\tEndHelix) + \hHelix });
      \coordinate (VerticalStart) at ({ \rx*cos(\tVertical) }, { \ry*sin(\tVertical) + \hVertical });

      \draw[very thick, orange, opacity=0.9] 
        (HelixEnd) .. controls +(5:0.5) and +(0,-0.4) .. (VerticalStart);

      \draw[very thick, orange, opacity=0.9] 
        (VerticalStart) -- ({ \rx*cos(\tVertical) }, { \h + \ry*sin(\tVertical) });

    \end{tikzpicture}
    }
  \end{minipage}
  \hfill
  \begin{minipage}{0.25\textwidth}
    \centering
    \resizebox{\linewidth}{!}%
    {%
      \begin{tikzpicture}

          \draw[thick, gray!80!black] (0,0) circle (3cm);
          
          \draw[thick, black!80] (0,0) circle (0.5cm);
          
          \coordinate (BlueTop) at (85:3cm);
          \coordinate (BlueBottom) at (-65:3cm);
          
          \coordinate (OrangeTop) at (65:3cm);
          \coordinate (OrangeBottom) at (-75:3cm); 
          
          \draw[very thick, cyan!80!black] (BlueTop) to[out=-70, in=100] (0.8,0.8) to[out=-80, in=95] (BlueBottom);
          \draw[very thick, orange!90!black] (OrangeTop) to[out=-100, in=75] (0.8,0.8) to[out=-105, in=100] (OrangeBottom);
          
          \fill[black!85] (-87:3cm) circle (4.5pt);
          
          \fill[black!85] (55:3cm) circle (4pt);

      \end{tikzpicture}
    }
  \end{minipage}

  \caption{Two curves on the cylinder. Note that when computing homs 
  from orange to blue we wrap the orange curve around the cylinder 
  until we hit a stop. Similarly, we are allowed to perform a Hamiltonian 
  isotopy (one which keeps the signed area between the curves equal to $0$).
  The last image give a simple example of this, where the orange curve is a 
  Hamiltonian deformation of the blue curve.}
\label{fig:Perturbation}
\end{figure}

We will compute the Fukaya category of a marked surface $S$ by fixing 
enough ``generating arcs'' to cut $S$ into topological disks.
We will be able to write every curve with a 
local system as an iterated extension of these generating arcs 
(equipped with trivial rank $1$ local systems). Combinatorially this will
correspond to homotoping a given curve to follow a sequence of generating arcs,
and algebraically it corresponds to the formation of an interesting 
dg-module\footnote{This is just a graded module $M$ equipped with a degree 
$1$ endomorphism $d : M \to M[1]$. Every chain complex $C^\bullet$ gives a 
dg-module by considering $\bigoplus_i C^i$ with differential $\bigoplus_i d^i$, 
but for more general graded algebras dg-modules are the correct thing to 
consider. See~\cite{kellerDifferentialGradedCategories2007} for an 
introduction.}.
Of course, there might be multiple choices for this sequence,
and this ambiguity is controlled by the presence of ``higher multiplication'' 
in the $A_\infty$-sense. If there is no such ambiguity then we say our surface 
is ``dissected'' and the higher $A_\infty$ multiplications all vanish. 
We'll first give a definition, 
then we'll give many examples which will hopefully help make this more clear.

\bigskip

All of our arcs and curves will be immersed, disjoint from $M$, and disjoint 
from the boundary except for endpoints of arcs, which must transversely end on 
the boundary.

\bigskip

\begin{defn}\label{def:GeneratingFamily}
  A \myimportant{Generating Family} for $(S,M)$ is a family of nonintersecting, 
  nonisotopic, simple arcs in $S$ so that every other such arc is homotopic 
  to a sequence of these.
\end{defn}

\begin{defn}\label{def:Dissection}
  A \myimportant{Dissection} for $(S,M)$ is a family of nonintersecting 
  nonisotopic simple arcs in $S$ which cuts the surface into 

  \begin{itemize}
    \item polygons, each of which has exactly one edge which is not
      in the family. This edge necessarily lies in the boundary.
    \item polygons where every edge is in the family, containing exactly 
      one puncture
  \end{itemize}
\end{defn}

\begin{figure}
  \centering
  
  \begin{minipage}{0.25\textwidth}
    \centering
    \hfill
    \resizebox{\linewidth}{!}%
    {%
      \begin{tikzpicture}

          \def\R{3}

          \draw[thick, black] (0,0) circle (\R);

          \draw[very thick, cyan, line cap=round] (180:\R) -- (0:\R);
          
          \draw[very thick, cyan, line cap=round] (160:\R) to[out=-25, in=-110] (100:\R);
          
          \draw[very thick, cyan, line cap=round] (195:\R) to[out=25, in=100] (260:\R);

          \node[fill=black, circle, minimum size=4.5mm] at (45:\R) {};
          \node[fill=black, circle, minimum size=4.5mm] at (140:\R) {};
          \node[fill=black, circle, minimum size=4.5mm] at (215:\R) {};
          \node[fill=black, circle, minimum size=4.5mm] at (315:\R) {};

      \end{tikzpicture}
    }
  \end{minipage}
  \hfill
  \begin{minipage}{0.25\textwidth}
    \centering
    \resizebox{\linewidth}{!}%
    {%
      \begin{tikzpicture}

          \newcommand{\rx}{4}
          \newcommand{\ry}{3}
          
          \draw[thick, black] (0,0) ellipse ({\rx} and {\ry});

          \coordinate (Lnode) at (-1.2, 0.3);
          \draw[thick, fill=white] (Lnode) circle (0.45);
          \fill[black] ($(Lnode)+(0.05,0.35)$) circle (0.15); 

          \coordinate (Rnode) at (1.2, 0.3);
          \draw[thick, fill=white] (Rnode) circle (0.45);
          \fill[black] ($(Rnode)+(-0.05,0.35)$) circle (0.15); 

          \begin{scope}[very thick, cyan]
              \draw (240:{\rx} and {\ry}) to[out=80, in=-80] ($(Lnode)+(0,-0.45)$);
              
              \draw ($(Lnode)+(-0.43, -0.15)$) to[out=170, in=-100] (145:{\rx} and {\ry});
              
              \draw ($(Lnode)+(0.4, -0.2)$) to[out=40, in=-110] (105:{\rx} and {\ry});
          \end{scope}

          \begin{scope}[very thick, cyan]
              \draw (265:{\rx} and {\ry}) to[out=95, in=-100] ($(Rnode)+(0,-0.45)$);
              
              \draw ($(Rnode)+(-0.4, -0.2)$) to[out=140, in=-70] (75:{\rx} and {\ry});
              
              \draw ($(Rnode)+(0.43, -0.15)$) to[out=10, in=-115] (40:{\rx} and {\ry});
          \end{scope}

          \fill[black!85] (90:{\rx} and {\ry}) circle (0.18);   
          \fill[black!85] (210:{\rx} and {\ry}) circle (0.18);  
          \fill[black!85] (320:{\rx} and {\ry}) circle (0.18);  

      \end{tikzpicture}
    }
  \end{minipage}
  \hfill
  \begin{minipage}{0.25\textwidth}
    \centering
    \resizebox{\linewidth}{!}%
    {%
      \begin{tikzpicture}[scale=1.5]
          \tikzset{
              blackpen/.style={very thick, draw=black!85, line cap=round, line join=round},
              arrowpen/.style={thick, draw=black!85, line cap=round, line join=round},
              bluepen/.style={very thick, color=cyan, line width=1.2pt, line cap=round, line join=round},
              orangepen/.style={very thick, color=orange, line width=1.2pt, line cap=round, line join=round}
          }

          \draw[blackpen] (0, 2) ellipse (1.6 and 0.3);
          
          \draw[blackpen] (-1.6, 2) .. controls (-1.8, 0.5) and (-1.2, -2) .. (0, -2.2) 
                          .. controls (1.2, -2) and (1.8, 0.5) .. (1.6, 2);

          \filldraw[black!85] (1.25, 1.8) circle (1.5pt);

          \draw[blackpen] (0.2, 0.3) .. controls (0.45, -0.1) and (0.45, -0.6) .. (0.1, -1.0);
          \draw[blackpen] (0.25, 0.1) .. controls (0.05, -0.2) and (0.05, -0.5) .. (0.25, -0.8);

          
          \draw[bluepen, dashed] (-1.41, -0.24) .. controls (-1.6, -0.18) and (0.3, -0.18) .. (0.12, -0.24);

          \draw[bluepen] (-1.1, 1.78) -- (-1.1, 0) arc (0:-105:0.25);

          \draw[bluepen] (-0.3, 1.70) -- (-0.3, -0.09) arc (180:285:0.2) -- (0.12, -0.24);

          \draw[bluepen] (-0.7, 1.73) -- (-0.7, -0.2)
              .. controls (-0.7, -1.7) and (0.9, -1.7) ..
              (0.9, -0.2) -- (0.9, 1.75);

      \end{tikzpicture}
    }
  \end{minipage}
  \caption{Some example dissections. When flowing along the boundary, we 
  always keep the surface on our right. So when computing homs between 
  these generating curves in the pair of pants we flow clockwise around the 
  outer boundary and counterclockwise around the inner boundary (until we 
  hit a stop)}
  \label{fig:Dissections}
\end{figure}

Finding generating families and dissections is usually not difficult to do by 
hand. Draw your surface in the plane (using a $2g$-gon for a surface of 
genus $g$, as usual), then add arcs until your surface is cut into disks.
Then add more edges to cut these 
disks into smaller disks each having exactly one boundary edge.
Figure~\ref{fig:BigDissection} shows a computation which is hopefully 
complicated enough to emphasize that with practice this is not so 
difficult\footnote{Throughout this paper we will use 
the convention (standard in the gentle algebra literature) that a dotted line
connecting two adjacent arrows of a quiver indicates that those arrows 
compose to $0$.}.

\begin{figure}
  \centering
  \begin{minipage}{0.4\textwidth}
    \centering
    \resizebox{\linewidth}{!}%
    {%
      \begin{tikzpicture}[scale=1.2, thick]

          \coordinate (vTR) at (22.5:4);  
          \coordinate (vT1) at (67.5:4);  
          \coordinate (vT2) at (112.5:4); 
          \coordinate (vTL) at (157.5:4); 
          \coordinate (vBL) at (202.5:4); 
          \coordinate (vB1) at (247.5:4); 
          \coordinate (vB2) at (292.5:4); 
          \coordinate (vBR) at (337.5:4); 

          \draw[color=cyan] (vT1) -- (vT2) -- (vTL) -- (vBL) -- (vB1) -- (vB2) -- (vBR) -- (vTR) -- cycle;

          \node[above] at (90:3.7) {\Large A};
          \node[above right] at (45:3.8) {\Large B};
          \node[right] at (0:3.8) {\Large A};
          \node[below right] at (315:3.8) {\Large D};
          \node[below] at (270:3.8) {\Large C};
          \node[below left] at (225:3.8) {\Large D};
          \node[left] at (180:3.8) {\Large C};
          \node[above left] at (135:3.8) {\Large B};

          \coordinate (c1) at (-1.5, 0.53); 
          \coordinate (c2) at (1.5, 1.2);
          \coordinate (c3) at (1.2, -1.2);
          
          \draw (c2) circle (0.4);
          \draw (c3) circle (0.4);

          \draw (vTL) -- ([shift=(90:1.0)]c1) arc (90:-90:1.0) to[out=180, in=-60] (vTL);

          \coordinate (p1) at ([shift=(110:1.06)]c1); 
          \coordinate (p2) at ([shift=(-100:1.0)]c1); 
          \coordinate (p3) at ([shift=(170:0.4)]c2);  
          \coordinate (p4) at ([shift=(170:0.4)]c3);  

          \fill (p1) circle (2.5pt);
          \fill (p2) circle (2.5pt);
          \fill (p3) circle (2.5pt);
          \fill (p4) circle (2.5pt);

    
          \draw[color=orange] ([shift=(85:1.0)]c1) to[out=65, in=105, looseness=1.6] (vBR);
          
          \draw[color=orange] ([shift=(20:1.0)]c1) to[out=20, in=220] ([shift=(220:0.4)]c2);
          \draw[color=orange] ([shift=(-40:0.4)]c2) to[out=-30, in=125, looseness=0.8] (vBR);

          \draw[color=orange] ([shift=(-15:1.0)]c1) to[out=-25, in=140] ([shift=(140:0.4)]c3);

          \draw[color=orange] ([shift=(-70:1.0)]c1) to[out=-70, in=215, looseness=1.2] (vBR);

          \node at (0, 2.5) {\Large E};
          \node at (-0.2, 1.2) {\Large F};
          \node at (2.4, -0.15) {\Large G};
          \node at (0.2, -0.6) {\Large H};
          \node at (-0.8, -1.6) {\Large I};

      \end{tikzpicture}
    }
  \end{minipage}
  \hfill
  \begin{minipage}{0.4\textwidth}
    \centering
    \[\begin{tikzcd}
      B & {} & C \\
      A && D \\
      E & G & I \\
      F && H
      \arrow[""{name=0, anchor=center, inner sep=0}, from=1-1, to=2-1]
      \arrow[""{name=1, anchor=center, inner sep=0}, shift right=3, from=1-3, to=2-3]
      \arrow[""{name=2, anchor=center, inner sep=0}, shift left=3, from=1-3, to=2-3]
      \arrow[""{name=3, anchor=center, inner sep=0}, shift right=3, from=2-1, to=1-1]
      \arrow[""{name=4, anchor=center, inner sep=0}, shift left=3, from=2-1, to=1-1]
      \arrow[""{name=5, anchor=center, inner sep=0}, from=2-3, to=1-3]
      \arrow[""{name=6, anchor=center, inner sep=0}, from=2-3, to=3-3]
      \arrow[""{name=7, anchor=center, inner sep=0}, from=3-1, to=2-1]
      \arrow[""{name=8, anchor=center, inner sep=0}, from=3-1, to=4-1]
      \arrow[""{name=9, anchor=center, inner sep=0}, from=3-2, to=3-1]
      \arrow[""{name=10, anchor=center, inner sep=0}, from=3-2, to=4-1]
      \arrow[""{name=11, anchor=center, inner sep=0}, from=3-3, to=3-2]
      \arrow[""{name=12, anchor=center, inner sep=0}, from=4-1, to=4-3]
      \arrow[""{name=13, anchor=center, inner sep=0}, from=4-3, to=3-3]
      \arrow[shift right=2, dotted, no head, from=0, to=4]
      \arrow[shift right=2, dotted, no head, from=2, to=5]
      \arrow[curve={height=-6pt}, dotted, no head, from=2, to=6]
      \arrow[shift left=2, dotted, no head, from=3, to=0]
      \arrow[shift left=2, dotted, no head, from=5, to=1]
      \arrow[curve={height=-6pt}, dotted, no head, from=7, to=4]
      \arrow[shift right=3, curve={height=-6pt}, between={0.2}{0.8}, dotted, no head, from=9, to=8]
      \arrow[shift right=5, curve={height=-6pt}, between={0.3}{1}, dotted, no head, from=10, to=12]
      \arrow[shift right=2, curve={height=-6pt}, between={0.2}{0.9}, dotted, no head, from=11, to=10]
      \arrow[shift right=2, curve={height=-6pt}, between={0.1}{0.9}, dotted, no head, from=13, to=11]
    \end{tikzcd}\]
  \end{minipage}

  \caption{A dissection of a genus $2$ surface with $3$ boundary components, 
  with two stops, one stop, one stop, respectively. To dissect a genus 
  $g$ surface with boundary we first draw the surface as a $4g$-gon in the 
  plane as usual (shown here in blue) with boundary components circles 
  (shown in black). One of the boundary components should touch the shared 
  point of all the blue curves (here the component with two stops). Then 
  one should add extra curves (shown in orange) to cut the surface into 
  polygons each of which contains exactly one boundary component with stop.
  In this example, writing $\bullet$ for a boundary component with a stop 
  those polygons are $BABAE\bullet$, $EG\bullet F$, $FGIH\bullet H$, and 
  $IDCDC\bullet$. At right is the $A_\infty$-category $Q$ 
  (read: quiver with relations) associated to this dissection, 
  so that the Fukaya category of this surface with stops is equivalent to 
  the derived category of $Q$-modules.} 
  \label{fig:BigDissection}
\end{figure}

\section{The $\mathbb{Z}/2$-graded Fukaya Category}

From a generating family for $(S,M)$, we build an $A_\infty$-category $Q$ 
(a fancy quiver) which presents the Fukaya category $\Fuk(S,M)$ 
as its derived category of modules.
In this category, we consider oriented curves, and changing the orientation 
corresponds to a homological shift. 
Of course, if we reverse the orientation twice we get back 
where we started, so that this gives a ``$\mathbb{Z}/2$-graded 
derived category''. These have been well studied after
Bridgeland's seminal paper~\cite{bridgelandQuantumGroupsHall2013} where they 
are used to relate quantum groups to Hall algebras.

The vertices of our quiver come from arcs in our family 
equipped with a choice of orientation and a rank $1$ local system.
Note that since our arcs are contractible, every local system is isomorphic to 
the trivial one, so there's no need to include it in our notation.
Arrows between these will be given by intersection points, more precisely by 
intersection points that arise after flowing our arcs keeping the surface on 
the right of the flow direction\footnote{Alternatively, which keep 
the intersection point/boundary on the \emph{left}}. There's a convenient 
shorthand for this in terms of ``angles'' connecting our curves, as shown in 
Figure~\ref{fig:ComposableArrows}. In the $\mathbb{Z}/2$-graded case, the 
arrows are placed in degree $1$ if the oriented arcs meet tip-to-tail and in 
degree $0$ otherwise. See Figure~\ref{fig:OrientedDegrees}. 

\begin{figure}
  \centering
  \begin{minipage}{0.30\textwidth}
    \centering
    \resizebox{\linewidth}{!}%
    {%
     \begin{tikzpicture}[
          >=Stealth, 
          force/.style={cyan, very thick, ->},
          dim/.style={black, thick, ->}
      ]

          \draw[thick] (-1,0) -- (3.8,0);

          \foreach \x in {-0.5, 1.0, ..., 3} {
              \draw[thick, gray!80] (\x, 0) -- (\x + 1.2, 1.0);
          }

          \coordinate (P1) at (0, 0);
          \coordinate (P2) at (1.5, 0);
          \coordinate (P3) at (3, 0);
          \coordinate (Dot) at (3.5, 0); 

          \draw[force] (0, -2) -- (P1);
          \draw[force] (1.5, -2) -- (P2);
          \draw[force] (P3) -- (3, -2); 

          \draw[dim] (0, -0.4) -- node[below] {$a$} (1.5, -0.4);
          \draw[dim] (1.5, -0.4) -- node[below] {$b$} (3, -0.4);

          \fill[black!85, draw=black] (Dot) circle (0.1);

      \end{tikzpicture}
    }
  \end{minipage}%
  \hfill
  \begin{minipage}{0.30\textwidth}
    \centering
    \resizebox{\linewidth}{!}%
    {%
      \begin{tikzpicture}[
          >=stealth, 
          wall/.style={thick, black!80}, 
          blueline/.style={thick, color=cyan!70!blue}
      ]

          \def\xmin{-1}
          \def\xmax{4}
          \def\ytop{4}
          \def\ybot{0}

          \draw[wall] (\xmin, \ytop) -- (\xmax, \ytop);
          \foreach \x in {-0.5, 0.5, ..., 3.5} {
              \draw[wall] (\x, \ytop) -- ++(0.4, 0.6);
          }

          \draw[wall] (\xmin, \ybot) -- (\xmax, \ybot);
          \foreach \x in {-0.5, 0.5, ..., 3.5} {
              \draw[wall] (\x, \ybot) -- ++(-0.4, -0.6);
          }

          \draw[blueline, <-] (0, \ytop) -- (0, \ytop-1.5);
          \draw[blueline, loosely dashed] (0, \ytop-1.5) -- (0, \ytop-2.5);

          \draw[blueline, <-] (3, \ytop) -- (3, \ybot);

          \draw[blueline, <-] (1.5, \ybot) -- (1.5, \ybot+1.5);
          \draw[blueline, loosely dashed] (1.5, \ybot+1.5) -- (1.5, \ybot+2.5);

          \draw[->, thick] (0, 3.4) -- node[below, font=\large] {$a$} (3, 3.4);

          \draw[->, thick] (3, 0.6) -- node[above, font=\large] {$b$} (1.5, 0.6);

      \end{tikzpicture}
    }
  \end{minipage}
  \hfill
  \begin{minipage}{0.30\textwidth}
    \centering
    \resizebox{\linewidth}{!}%
    {%
      \begin{tikzpicture}[
          >={Stealth[length=3mm, width=2.5mm]},
          stream/.style={color=cyan!70!blue, very thick},
          boundary/.style={color=gray!80!black, thick},
          pencil/.style={color=gray!80!black}
      ]

      \draw[boundary] (-2.5, 4) -- (6.8, 4);

      \foreach \x in {-2, -1, 0, 1, 2, 3, 4, 5, 6} {
          \draw[boundary] (\x, 4) -- (\x+0.6, 4.8);
      }

      \node[circle, fill=black!85, inner sep=2pt] at (6.8, 4) {};


      \draw[stream, <-] (2, -1.5) -- (2, 4);

      \draw[stream] (-1.5, -2) -- (-1.5, 0) to[out=90, in=180] (0.5, 2.5) -- (4, 2.5);

      \draw[stream] (-0.5, -2) -- (-0.5, -0.5) to[out=90, in=180] (1, 1.5) -- (2.5, 1.5) to[out=0, in=-90] (4, 2.5);

      \draw[stream, ->] (4, 2.5) -- (4, 4);

      \draw[stream, ->, rounded corners=12pt] (4, 2.5) -- (6, 2.5) -- (6, 4);


      \node[circle, fill=black!85, inner sep=2.5pt] at (4, 2.5) {};
      \node[pencil, below right=1pt, font=\Large] at (4, 2.5) {$a$};

      \node[circle, fill=black!85, inner sep=2.5pt] at (2, 1.5) {};
      \node[pencil, below right=1pt, font=\Large] at (2, 1.5) {$b$};

      \end{tikzpicture}
    }
  \end{minipage}%
  \caption{There's a very convenient combinatorial 
  shorthand for maps between our generating arcs, which only intersect
  along the boundary after wrapping. These are given by paths (also called 
  ``angles'') between our arcs keeping the surface on the right.
  The homological degree is computed from the orientation 
  data as usual (so here $a$ has degree $0$ and $b$ has degree $1$), 
  and these compose when $a$ and $b$ are adjacent as arcs. So in the first 
  subfigure $a$ and $b$ are composable, while in the second subfigure they 
  are not (formally, the composite is defined to be $0$). 
  In the symplectic/geometric picture, these angles correspond to intersection
  points that arise after wrapping, and the composition comes from the presence
  of a holomorphic triangle as shown in the last subfigure.}
  \label{fig:ComposableArrows}
\end{figure}


\begin{figure}
  \centering

  \begin{minipage}{0.1\textwidth}
    \centering
    \resizebox{\linewidth}{!}%
    {%
      \begin{tikzpicture}

          \draw[orange, line width=2.5pt, -{Triangle[length=8pt, width=8pt]}] (0.6, -1.1) -- (-0.6, 1.1);
          
          \draw[cyan, line width=2.5pt, -{Triangle[length=8pt, width=8pt]}] (-0.6, -1.1) -- (0.6, 1.1);
          
          \fill[black] (0, 0) circle (4.5pt);
          
          \draw[black, line width=1.2pt, ->] (0.26, -0.48) arc (-61.4:61.4:0.55);
          \draw[black, line width=1.2pt, ->] (-0.26, 0.48) arc (118.6:241.4:0.55);
          
          \node[black, font=\large\bfseries] at (0.8, 0) {$1$};

      \end{tikzpicture}
    }
  \end{minipage}%
  \quad\quad\quad\quad
  \begin{minipage}{0.1\textwidth}
    \centering
    \resizebox{\linewidth}{!}%
    {%
      \begin{tikzpicture}

          \draw[orange, line width=2.5pt, -{Triangle[length=8pt, width=8pt]}] (-0.6, 1.1) -- (0.6, -1.1);
          
          \draw[cyan, line width=2.5pt, -{Triangle[length=8pt, width=8pt]}] (-0.6, -1.1) -- (0.6, 1.1);
          
          \fill[black] (0, 0) circle (4.5pt);
          

          \draw[black, line width=1.2pt, ->] (0.26, -0.48) arc (-61.4:61.4:0.55);
          \draw[black, line width=1.2pt, ->] (-0.26, 0.48) arc (118.6:241.4:0.55);
          
          \node[black, font=\large\bfseries] at (0.7, 0) {$0$};

      \end{tikzpicture}
    }
  \end{minipage}%
  \hfill

  \caption{We compute the degree of an intersection point based on the 
  relative orientations of the intersecting curves.}
  \label{fig:OrientedDegrees}
\end{figure}

Composition comes from putting two 
adjacent angles together into a bigger angle (and is $0$ if this is not possible), 
and the higher $A_\infty$-operations come from closed polygons. 
See~\cite[Chapter 9.1]{bocklandtGentleIntroductionHomological2021} for more 
information.

For example, the Fukaya category of a disk with $n$ marked boundary points is the 
$A_{n-1}$-quiver by generalizing Figure~\ref{fig:AnQuiver}

\begin{figure}
  \centering
  \begin{minipage}{0.3\textwidth}
    \centering
    \resizebox{\linewidth}{!}%
    {%
      \begin{tikzpicture}[>=Stealth]

          \draw[thick, black] (0,0) circle (3cm);

          \fill[black] (135:3) circle (4pt);
          \fill[black] (45:3) circle (4pt);
          \fill[black] (-45:3) circle (4pt);
          \fill[black] (-135:3) circle (4pt);

          \draw[thick, cyan, ->, name path=lineB] (180:3) -- (0:3) 
              node[right, text=cyan, scale=1.5] {B};
          
          \draw[thick, cyan, ->, name path=curveA] (200:3) to[bend left=35] (-90:3) 
              node[below, text=cyan, scale=1.5] {A};
          
          \draw[thick, cyan, ->, name path=curveC] (160:3) to[bend right=35] (90:3) 
              node[above, text=cyan, scale=1.5] {C};
          
          \draw[thick, orange, ->] (-75:3) to[bend left=35] (-15:3) 
              node[right, text=orange, scale=1.5] {X};
          
          \path[name path=vert] (-2.2, -2) -- (-2.2, 2);
          
          \path[name intersections={of=curveA and vert, by=A_pt}];
          \path[name intersections={of=lineB and vert, by=B_ptA}];
          \draw[thick, ->] (A_pt) -- (B_ptA) 
              node[midway, right, scale=1.3] {a};
          
          \path[name intersections={of=lineB and vert, by=B_ptB}];
          \path[name intersections={of=curveC and vert, by=C_pt}];
          \draw[thick, ->] (B_ptB) -- (C_pt) 
              node[midway, right, scale=1.3] {b};
      \end{tikzpicture}
    }
  \end{minipage}
  \hfill
  \begin{minipage}{0.3\textwidth}
    \centering
    \resizebox{\linewidth}{!}%
    {%
      \begin{tikzpicture}
          \draw[thick] (0,0) circle (3);

          \draw[cyan, ultra thick, -Stealth] (-3,0) -- (3,0) node[right, font=\Large\bfseries] {B};
          \draw[cyan, ultra thick, -Stealth] (0,-3) -- (0,3) node[above, font=\Large\bfseries] {Y};

          \fill (0,0) circle (3.5pt);
          \node[below left, font=\large] at (-0.05,-0.05) {$p$};

          \draw[thick, ->] (100:0.4) arc (100:170:0.4);
          \draw[thick, ->] (280:0.4) arc (280:350:0.4);

          \foreach \angle in {45, 135, 225, 315} {
              \fill (\angle:3) circle (3.5pt);
          }

          \draw[orange, ultra thick, -Stealth] (-2.95,0.52) arc (270:360:2.5);
          
          \draw[orange, ultra thick, -Stealth] (0.52,-2.95) arc (180:90:2.5);

          \node[orange, above left, font=\Large\bfseries] at (-3,0.3) {C};
          \node[orange, below, font=\Large\bfseries] at (0.5,-3.2) {X};
      \end{tikzpicture}
    }
  \end{minipage}
  \caption{A dissection (shown in blue) of the disk with $4$ stops, showing that its 
  Fukaya category is equivalent to representations of the $A_3$ quiver
  $C \overset{b}{\to} B \overset{a}{\to} A$. Note that the quiver arrows 
  face the opposite direction from the maps of curves, essentially because 
  the Yoneda embedding is contravariant and our generating curves become the 
  projective representations of our quiver. The generator 
  orientations were chosen to put the entire algebra in degree $0$.
  The orange curve $X$ is homotopic to $A$, with the opposite orientation, 
  glued to $B$ along the angle $a$. This means that algebraically $X$ is 
  represented by $\text{Cone}(a)$, that is, by the dg-module 
  $\left ( A[-1] \oplus B, \begin{pmatrix} 0 & 0 \\ a & 0 \end{pmatrix} \right )$.
  The second figure shows how $\text{Cone}(p)$ is computed as a resolution 
  of the intersection point.}
  \label{fig:AnQuiver}
\end{figure}

The cone of a degree $1$ angle comes from gluing its two arcs together into 
a longer curve and this tells us how to compute the algebraic representative of a 
general arc in $(S,M)$ -- An oriented arc in our generating family is represented by 
the projective module associated to that arc's vertex in our quiver if the 
orientations agree, or the shift of that projective if the orientations disagree.
For an arc outside our generating family, 
we isotope it so that it lies along our generating family,
orienting the arcs consistently, then repeatedly take cones along the 
degree $1$ maps between these. A simple example is 
shown in figure~\ref{fig:AnQuiver}, where the curves $A$, $B$, and $C$ 
are represented by the projectives at the corresponding vertices\footnote{For 
certain readers, it might help to recognize these projective modules as 
the images of the Yoneda embedding $Q^\text{op} \hookrightarrow Q\text{-mod}$}. 
These are, respectively, the representations $0 \to 0 \to k$, $0 \to k \to k$, 
and $k \to k \to k$. Then since $Q$ is concentrated in degree $0$ dg-modules 
are just chain complexes, and the orange curve $X$ can be represented by 
the complex of $Q$-modules $(0 \to 0 \to k) \to (0 \to k \to k)$ where the 
differential is the obvious inclusion. Of course, this complex is quasi-equivalent 
to the module $0 \to k \to 0$. This is not an accident, since in general 
the curves on our surface (with local systems) will be in bijection with the 
indecomposables of (the derived category of) $Q$-modules. It makes a great 
exercise to compute the algebraic representatives of the $6$ curves on the 
disk with $4$ stops and check that they really are the $6$ indecomposables for
the $A_3$ quiver. In general, the cone of an interior marked point is 
given by the resolution of that marked point, and another great exercise is
to verify the relation shown in the right subfigure of 
Figure~\ref{fig:AnQuiver}. It says that $\text{Ext}^\bullet(Y,B)$ should have 
a single generator $p$ in degree $1$, and that $\text{Cone}(p)$ should be 
$X \oplus C$. For those who want to relate the combinatorics to the 
symplectic geometry, note that our arcs become intersection points 
after wrapping, and resolving the intersection point gives the glued curve 
plus a curve which we can isotope away to infinity.

The situation for a closed curve is similar, but now that it's noncontractible 
we're forced to keep track of a specific choice of local system. This is the 
choice of a finite dimensional $k[X^\pm]$-module, where we think of the action 
of $X$ as the \emph{monodromy} action when we go once around our closed loop.
As for arcs, this local system will decompose into a finite direct sum of 
indecomposable $k[X^\pm]$-modules, which (after choosing a basis) are in bijection 
with invertible companion matrices~\cite[Chapter VI.7]{aluffiAlgebraChapter02009}.

To build a dg-module from a closed curve with monodromy $M \in \text{GL}_r(k)$, 
we proceed as in the arc case by deforming our curve to lie along 
the generating curves and using the angles between generating curves as our 
differential. The difference is that now we have an extra angle $\alpha_*$
connecting our last generating arc back to our first generating arc. We 
take the direct sum of $r$ many copies of our arcs, using 
$\text{Id}_r \otimes \alpha$ for the differentials along the way, except 
for the final angle connecting our starting and ending curves, where we use 
$M \otimes \alpha_*$ instead.

  \begin{figure}
    \centering
    \begin{minipage}{0.3\textwidth}
      \centering
      \resizebox{\linewidth}{!}%
      {%
        \begin{tikzpicture}
          \tikzset{%
              grey line/.style={draw=black, line width=1.5pt},
              orange line/.style={draw=orange, line width=2pt},
              blue line/.style={draw=cyan, line width=2pt, -{Stealth[length=12pt, width=10pt, bend]}},
              label arrow/.style={draw=black!85, line width=1pt, -{Stealth[length=6pt, width=5pt]}}
          }

          \draw[grey line] (0,0) circle (4.5cm);
          \fill[black] (0, -4.5) circle (0.22cm);

          \draw[grey line] (0,0) circle (1.0cm);
          \fill[black!85] (0, 1.0) circle (0.18cm);

          \draw[orange line] (0,0) circle (2.75cm);
          \draw[orange line, -{Stealth[length=12pt, width=10pt]}] (0.1, 2.75) -- (-0.1, 2.75);

          \draw[blue line] (-0.8, 4.4) .. controls (-1.1, 0.5) and (-1.9, -3.8) .. (240:1.0);
          
          \draw[blue line] (0.8, 4.4) .. controls (1.1, 0.5) and (1.9, -3.8) .. (300:1.0);

          \draw[label arrow, ->] (-0.8, 3.8) to[bend left=15] node[below=1pt] {\Huge $a$} (0.8, 3.8);

          \draw[label arrow] (240:1.8) to[bend right=25] node[below=2pt] {\Huge $b$} (300:1.8);
      \end{tikzpicture}
    }
  \end{minipage}

    \caption{Two generating arcs ($L\text{eft}$ and $R\text{ight}$) are 
    shown in blue with the two maps $a,b : L \to R$ so that we see the 
    Fukaya category of this surface is the derived category of representations 
    of the Kronecker quiver $\bullet \rightrightarrows \bullet$. In this 
    section we mean the $\mathbb{Z}/2$-graded derived category, but of course 
    this will be true in the $\mathbb{Z}$-graded case too, given an appropriate 
    line field.}
  \label{fig:Monodromy}
\end{figure}

For example, consider the closed curve shown in orange in Figure~\ref{fig:Monodromy}, 
with monodromy given by the $3 \times 3$ companion matrix 
$\begin{pmatrix} 0 & 0 & C_0 \\ 1 & 0 & C_1 \\ 0 & 1 & C_2 \end{pmatrix}$.
We deform this curve to lie along the generating 
arcs $L$ and $R$ (with the opposite orientation, corresponding to a 
homological shift). Starting at the top and moving counterclockwise, we first 
pass the angle $b$, and then we close our loop with the angle $a$. 
Since our local system is $3$-dimensional, we consider the dg-module  
$L^{\oplus 3} \oplus {R[1]}^{\oplus 3}$ with differential 
$\text{Id}_3 \otimes b + M \otimes a$. Concretely this is the matrix

\[ 
\begin{pNiceArray}{ccc|ccc}
  \Block{3-3}<\Large>{\mathbf{0}}
  & & & \Block{3-3}<\Large>{\mathbf{0}} \\
  & & & & &\\
  & & & & &\\
  \hline
  b & 0 & C_0 a & \Block{3-3}<\Large>{\mathbf{0}} \\
  a & b & C_1 a & & & \\
  0 & a & C_2 a + b & & &
\end{pNiceArray}  
\]

In the interest of space, we will abbreviate such matrices as a sum of 
their entries. For example, the above matrix will be written 
$b^{41} + b^{52} + b^{63} + a^{51} + a^{62} + 
C_0 a^{43} + C_1 a^{53} + C_2 a^{63}$, or 
(with reference to the $2 \times 2$ block diagonal structure)
as $\text{Id}_3 \otimes b^{21} + M \otimes a^{21}$.

\begin{figure}
  \centering
  \begin{minipage}{0.4\textwidth}
    \centering
    \resizebox{\linewidth}{!}%
    {%
      \begin{tikzpicture}[
          thick,
          >={Stealth[length=2.5mm, width=2mm]},
          bluepen/.style={color=cyan, thick},
          orangepen/.style={color=orange, thick},
          graypen/.style={color=black!85, thick}
      ]

      \draw[graypen] (0,0) ellipse (4.5cm and 2.5cm);
      \filldraw[graypen] (0, 2.5) circle (2.5pt); 

      \draw[orangepen, decoration={
          markings,
          mark=at position 0.32 with {\arrow[scale=1.2]{>}}
      }, postaction={decorate}] (0,0) ellipse (3.4cm and 1.1cm);

      \draw[bluepen, <-] (1.8, -2.29) -- (1.8, 0.5) arc(180:0:0.6 and 1.0) -- (3.0, -1.86);
      \node at (2.4, 0) {\Large $\times$};
      \node[bluepen, font=\Large] at (2.4, -3.2) {C};

      \draw[bluepen, <-] (-0.6, -2.48) -- (-0.6, 0.5) arc(180:0:0.6 and 1.0) -- (0.6, -2.48);
      \node at (0, 0) {\Large $\times$};
      \node[bluepen, font=\Large] at (0, -3.2) {B};

      \draw[bluepen, <-] (-3.0, -1.86) -- (-3.0, 0.5) arc(180:0:0.6 and 1.0) -- (-1.8, -2.29);
      \node at (-2.4, 0) {\Large $\times$};
      \node[bluepen, font=\Large] at (-2.4, -3.2) {A};

      \def\arrY{-1.5}
      \draw[graypen, ->] (3.0, \arrY) -- (1.8, \arrY) node[midway, above=1pt, text=black!80] {$a$};
      \draw[graypen, ->] (1.8, \arrY) -- (0.6, \arrY) node[midway, above=1pt, text=black!80] {$b$};
      \draw[graypen, ->] (0.6, \arrY) -- (-0.6, \arrY) node[midway, above=1pt, text=black!80] {$c$};
      \draw[graypen, ->] (-0.6, \arrY) -- (-1.8, \arrY) node[midway, above=1pt, text=black!80] {$d$};
      \draw[graypen, ->] (-1.8, \arrY) -- (-3.0, \arrY) node[midway, above=1pt, text=black!80] {$e$};

      \end{tikzpicture}
    }
  \end{minipage}
  \hfill
  \begin{minipage}{0.5\textwidth}
    \centering
    \resizebox{\linewidth}{!}%
    {%
      \begin{tikzpicture}[>=Stealth, auto, node distance=3cm, thick]

        \node (A) at (0,0) {\large $A$};
        \node (B) at (3,0) {\large $B$};
        \node (C) at (6,0) {\large $C$};

        \draw[->] (A) -- node[above] {$d$} (B);
        \draw[->] (B) -- node[above] {$b$} (C);

        \path[->] (C) edge[loop above, distance=1.8cm, out=125, in=55] node[above=0.1cm] {$a$} (C);
        \draw[dashed, thin] ($(C.center)+(-0.35cm,0.6cm)$) to[bend left=15] ($(C.center)+(0.35cm,0.6cm)$);

        \path[->] (B) edge[loop above, distance=1.8cm, out=125, in=55] node[above=0.1cm] {$c$} (B);
        \draw[dashed, thin] ($(B.center)+(-0.35cm,0.6cm)$) to[bend left=15] ($(B.center)+(0.35cm,0.6cm)$);

        \path[->] (A) edge[loop above, distance=1.8cm, out=125, in=55] node[above=0.1cm] {$e$} (A);
        \draw[dashed, thin] ($(A.center)+(-0.35cm,0.6cm)$) to[bend left=15] ($(A.center)+(0.35cm,0.6cm)$);

      \end{tikzpicture}
    }
  \end{minipage}
  \caption{All these maps live in degree $1$, and so we can take the cone 
  along any of them. Working in the $\mathbb{Z}/2$-graded setting, we can 
  also take the cone along an odd combination of them, like $edcba$, whereas 
  in the $\mathbb{Z}$-graded world this map would have degree $5$. This 
  corresponds to the fact that the orange curve admits a consistent 
  orientation, but will not admit a consistent grading against a line field
  that puts all these maps in degree $1$.}
\label{fig:ThreeLegPairOfPants}
\end{figure}

Similarly, say we have the curve in Figure~\ref{fig:ThreeLegPairOfPants}
with monodromy given by some $2 \times 2$ companion matrix $M$. Then we have 
a dg-module $A^{\oplus 2} \oplus B^{\oplus 2} \oplus C^{\oplus 2}$ 
with differential $\text{Id}_2 \otimes b^{21} + \text{Id}_2 \otimes d^{32}
+ M \otimes edcba^{31}$. Note that the map $edcba$ is degree $1$ (as 
needed to be part of a differential) since homological degrees are 
computed mod $2$.

You can compute the $\text{Ext}^\bullet$ groups between two curves $\alpha$ 
and $\beta$ by looking at their intersection points. Each interior intersection 
point gives one map in $\text{Ext}^\bullet(\alpha,\beta)$ and one map in 
$\text{Ext}^\bullet(\beta,\alpha)$, corresponding to the complementary
angles.
The cone of an internal degree $1$ intersection is given by the direct sum 
of the two curves in the resolution of this intersection, see for instance,
Figure~\ref{fig:AnQuiver} and Figure~\ref{fig:PuncturedTorusExample}.
Lastly, we must remember to handle the ``wrapping'' for a 
cylindrical boundary without stops, as in Figure~\ref{fig:PantsVariations}. 

\section{Line Fields and the $\mathbb{Z}$-graded Fukaya Category}
In this section we'll explain how to equip $(S,M)$ with an extra piece 
of data (a \emph{line field}) which will let us put a $\mathbb{Z}$-grading 
on $\Fuk(S,M)$ and thus recover the usual notion of derived category. 
Instead of equipping our curves with an \emph{orientation} (which has a 
$\mathbb{Z}/2$-torsor worth of options), we equip them with a \emph{grading} 
(which has a $\mathbb{Z}$-torsor worth of options). Then our objects will be 
graded arcs, which will let us put intersections points into $\mathbb{Z}$-graded
degrees as well.

Officially, a line field $\eta$ in a surface $S$ is a section of the projectivized 
tangent bundle $\mathbb{P}TS$, but it's often easier to think in terms of 
the corresponding foliation. This is a decomposition of $S$ into a 
disjoint union of curves, as in Figure~\ref{fig:LineField}

\begin{figure}
  \centering
  \begin{minipage}{0.5\textwidth}
    \centering
    \resizebox{\linewidth}{!}%
    {%
      \begin{tikzpicture}

          \def\Rout{3.6}
          \def\Rmid{2.2}
          \def\Rin{1.1}

          \begin{scope}
              \path[clip] (0:\Rout) arc (0:360:\Rout) -- (0:\Rin) arc (360:0:\Rin) -- cycle;
              
              \foreach \y in {-3.3,-3.0,...,3.6} {
                  \draw[gray!50, thick] (-\Rout,\y) -- (\Rout,\y);
              }
          \end{scope}

          \draw[black, line width=2pt] (0,0) circle (\Rout);
          \draw[black, line width=2pt] (0,0) circle (\Rin);

          \draw[orange, line width=2pt] (0,0) circle (\Rmid);

          \foreach \ang in {-22, 20, 60, 101, 141, 182, 230, 290} {
              \draw[black, thick, <-, >=Stealth] 
                  (\ang+20:\Rmid) ++(\ang+120:0.35) 
                  arc (\ang+120:0:0.25);
          }

          \draw[black, thick, <-, >=Stealth] (2.067, -.752) 
              to[out=15, in=180] (4.8, 1.2) 
              node[right, font=\sffamily\bfseries] {\Large uh oh!};
      \end{tikzpicture}
    }
  \end{minipage}
  \hfill
  \begin{minipage}{0.35\textwidth}
    \centering
    \resizebox{\linewidth}{!}%
    {%
      \begin{tikzpicture}

          \def\rInner{1.2}
          \def\rMiddle{2.2}
          \def\rOuter{3.4}

          \foreach \angle in {0,30,...,330} {
              \draw[thick, gray] (\angle:\rInner) -- (\angle:\rOuter);
          }

          \draw[ultra thick, orange] (0,0) circle (\rMiddle);

          \foreach \angle in {0,30,...,330} {
              \draw[thick, -{Stealth[length=6pt, width=4pt]}] (\angle:\rMiddle+0.5) arc (\angle:\angle+90:0.5);
          }

          \draw[ultra thick, black] (0,0) circle (\rInner);
          \draw[ultra thick, black] (0,0) circle (\rOuter);

      \end{tikzpicture}
    }
  \end{minipage}

  \caption{Two line fields on the annulus. In the first, the closed curve 
  does not admit a consistent grading since when we close the loop we differ
  from where we started by two half turns. In the second line field, however,
  the closed curve admits a consistent grading as shown.}
\label{fig:LineField}
\end{figure}

Then a \myimportant{grading} on a curve $\gamma$ is a continuous choice of 
homotopy at each point from the line field to the tangent line $\dot{\gamma}$. 
In the $\mathbb{Z}/2$-graded derived category, a homological shift was 
represented by a change in orientation. Now in the $\mathbb{Z}$-graded case a 
homological shift comes from adding or removing a half counterclockwise turn 
from the grading.

Note also that while every closed curve admits an \emph{orientation}, it's 
not the case that every closed curve admits a grading!
This means that (depending on our choice of line field) not every 
object in the $\mathbb{Z}/2$-graded Fukaya category lifts to an object in the 
$\mathbb{Z}$-graded Fukaya category, as shown in Figure~\ref{fig:LineField}.

As before, an internal intersection point gives an element of both 
$\text{Ext}^\bullet(\alpha, \beta)$ and $\text{Ext}^\bullet(\beta, \alpha)$,
though determining the grading is now slightly more involved. To determine 
the degree of the intersection point (viewed as a map $\alpha \to \beta$) 
we start at the line field, follow our homotopy to $\dot{\alpha}$, take the shortest 
counterclockwise path to $\dot{\beta}$, and then follow the \emph{inverse} 
homotopy from $\dot{\beta}$ back to the line field. The degree of our map is the 
signed number of counterclockwise half turns we made at the end of this process. 
As a helpful check, notice that the sum of the two degrees must be $1$.

For example, in Figure~\ref{fig:GradedIntersection} we compute the 
degree of the intersection point 
viewed as a map from the orange curve to the blue curve. To do this, we 
first follow the grading from the line field to the orange curve, then we 
take the shortest counterclockwise path to the blue curve, then undo the 
blue grading to end up back at the line field. We count that from start 
to end we've performed $1$ half turn, so that this map is degree $1$.
It's a nice exercise to compute the degree of this same intersection point,
viewed as a map from blue to orange, lies in degree $0$. This is consistent 
with the helpful check from earlier, since $1+0 = 1$.

\begin{figure}
  \centering
  \begin{minipage}{0.2\textwidth}
    \centering
    \resizebox{\linewidth}{!}%
    {%
      \begin{tikzpicture}

        \def\R{3}

        \draw[very thick] (0,0) circle (\R);

        \begin{scope}
          \clip (0,0) circle (\R);
          \foreach \y in {-2.5, -2.0, -1.5, -1.0, -0.5, 0, 0.5, 1.0, 1.5, 2.0, 2.5} {
            \draw[gray, thick] (-\R, \y) -- (\R, \y);
          }
        \end{scope}

        \draw[orange, very thick] (225:\R) -- (45:\R);

        \draw[cyan, very thick] (135:\R) -- (-45:\R);

        \draw[-stealth, thick] (-1.5, 1.5) + (0:0.3) arc (0:135:0.3);

        \draw[-stealth, thick] (1.5, 1.5) + (0:0.3) arc (0:225:0.3);

      \end{tikzpicture}
    }
  \end{minipage}
  \hfill
  \begin{minipage}{0.2\textwidth}
    \centering
    \resizebox{\linewidth}{!}%
    {%
      \begin{tikzpicture}

        \def\R{3}

        \draw[very thick] (0,0) circle (\R);

        \begin{scope}
          \clip (0,0) circle (\R);
          \foreach \y in {-2.5, -2.0, -1.5, -1.0, -0.5, 0, 0.5, 1.0, 1.5, 2.0, 2.5} {
            \draw[gray, thick] (-\R, \y) -- (\R, \y);
          }
        \end{scope}

        \draw[orange, very thick] (225:\R) -- (45:\R);

        \draw[cyan, very thick] (135:\R) -- (-45:\R);

        \draw[-stealth, thick] (-1.5, 1.5) + (0:0.3) arc (0:135:0.3);

        \draw[-stealth, thick] (1.5, 1.5) + (0:0.3) arc (0:225:0.3);

        \draw[-stealth, thick] (0, 0) + (0:0.3) arc (0:225:0.3);

      \end{tikzpicture}
    }
  \end{minipage}
  \hfill
  \begin{minipage}{0.2\textwidth}
    \centering
    \resizebox{\linewidth}{!}%
    {%
      \begin{tikzpicture}

        \def\R{3}

        \draw[very thick] (0,0) circle (\R);

        \begin{scope}
          \clip (0,0) circle (\R);
          \foreach \y in {-2.5, -2.0, -1.5, -1.0, -0.5, 0, 0.5, 1.0, 1.5, 2.0, 2.5} {
            \draw[gray, thick] (-\R, \y) -- (\R, \y);
          }
        \end{scope}

        \draw[orange, very thick] (225:\R) -- (45:\R);

        \draw[cyan, very thick] (135:\R) -- (-45:\R);

        \draw[-stealth, thick] (-1.5, 1.5) + (0:0.3) arc (0:135:0.3);

        \draw[-stealth, thick] (1.5, 1.5) + (0:0.3) arc (0:225:0.3);

        \draw[-stealth, thick] (0, 0) + (0:0.3) arc (0:225:0.3);
        \draw[-stealth, thick] (0, 0) + (225:0.45) arc (225:315:0.45);

      \end{tikzpicture}
    }
  \end{minipage}
  \hfill
  \begin{minipage}{0.2\textwidth}
    \centering
    \resizebox{\linewidth}{!}%
    {%
      \begin{tikzpicture}

        \def\R{3}

        \draw[very thick] (0,0) circle (\R);

        \begin{scope}
          \clip (0,0) circle (\R);
          \foreach \y in {-2.5, -2.0, -1.5, -1.0, -0.5, 0, 0.5, 1.0, 1.5, 2.0, 2.5} {
            \draw[gray, thick] (-\R, \y) -- (\R, \y);
          }
        \end{scope}

        \draw[orange, very thick] (225:\R) -- (45:\R);

        \draw[cyan, very thick] (135:\R) -- (-45:\R);

        \draw[-stealth, thick] (-1.5, 1.5) + (0:0.3) arc (0:135:0.3);

        \draw[-stealth, thick] (1.5, 1.5) + (0:0.3) arc (0:225:0.3);

        \draw[-stealth, thick] (0, 0) + (0:0.3) arc (0:225:0.3);
        \draw[-stealth, thick] (0, 0) + (225:0.45) arc (225:315:0.45);
        \draw[-stealth, thick] (0, 0) + (315:0.6) arc (315:180:0.6);

      \end{tikzpicture}
    }
  \end{minipage}
  \hfill
  \caption{Two graded curves near an intersection point. We compute the 
  degree of the intersection point (viewed as a hom from orange to blue) by 
  going from the line field to the orange curve (by the orange grading), 
  then to the blue curve (by the shortest counterclockwise path), then 
  back to the line field (by the blue grading backwards). The degree of the 
  intersection point is then the number of counterclockwise half turns, 
  so that the intersection point shown is placed in degree $1$.}
\label{fig:GradedIntersection}
\end{figure}

These kinds of computations with winding numbers are not difficult, but 
can be error prone in practice. Thankfully there is always an \emph{canonical} 
choice of line field associated to a dissection, 
described in~\cite{hausdorffcenterformathematicsGeometricModelBounded2020}
and~\cite[Section 2]{opperGeometricModelDerived2018},
where the algebra is placed entirely in degree $0$ and the degree computations 
are particularly combinatorial. In fact, more generally than degree $0$, 
for any grading on a gentle algebra, one can build a line field instantiating 
that grading~\cite[Theorem 3.11]{lekiliDerivedEquivalencesGentle2020}. 

\section{Examples}

\subsection{The Punctured Torus}

\begin{figure}
  \centering
  \begin{minipage}{0.24\textwidth}
    \centering
    \resizebox{\linewidth}{!}%
    {%
      \begin{tikzpicture}[scale=1.5]
          \tikzset{
              blackpen/.style={very thick, draw=black!85, line cap=round, line join=round},
              arrowpen/.style={thick, draw=black!85, line cap=round, line join=round},
              bluepen/.style={very thick, color=cyan, line width=1.2pt, line cap=round, line join=round},
              orangepen/.style={very thick, color=orange, line width=1.2pt, line cap=round, line join=round}
          }

          \draw[blackpen] (0, 2) ellipse (1.6 and 0.3);
          
          \draw[blackpen] (-1.6, 2) .. controls (-1.8, 0.5) and (-1.2, -2) .. (0, -2.2) 
                          .. controls (1.2, -2) and (1.8, 0.5) .. (1.6, 2);

          \filldraw[black!85] (1.25, 1.8) circle (1.5pt);

          \draw[blackpen] (0.2, 0.3) .. controls (0.45, -0.1) and (0.45, -0.6) .. (0.1, -1.0);
          \draw[blackpen] (0.25, 0.1) .. controls (0.05, -0.2) and (0.05, -0.5) .. (0.25, -0.8);

          
          \draw[bluepen, dashed] (-1.41, -0.24) .. controls (-1.6, -0.18) and (0.3, -0.18) .. (0.12, -0.24);

          \draw[bluepen] (-1.1, 1.78) node[above] {A} -- (-1.1, 0) arc (0:-105:0.25);

          \draw[bluepen, <-] (-0.3, 1.70) -- (-0.3, -0.09) arc (180:285:0.2) -- (0.12, -0.24);

          \draw[orangepen, <-] (-0.7, 1.73) node[above] {B} -- (-0.7, -0.2)
              .. controls (-0.7, -1.7) and (0.9, -1.7) ..
              (0.9, -0.2) -- (0.9, 1.75);


          \draw[arrowpen, ->] (-1.1, 1.2) -- (-0.7, 1.2) node[above, midway] {a};
          \draw[arrowpen, ->] (-0.7, 1.2) -- (-0.3, 1.2) node[above, midway] {b};
          \draw[arrowpen, ->] (-0.3, 1.2) -- (0.9, 1.2) node[above, midway] {c};

      \end{tikzpicture}
    }
  \end{minipage}
  \hfill
  \begin{minipage}{0.24\textwidth}
    \centering
    \[\begin{tikzcd}
        A && B
        \arrow[""{name=0, anchor=center, inner sep=0}, "a", shift left, curve={height=-12pt}, from=1-1, to=1-3]
        \arrow[""{name=1, anchor=center, inner sep=0}, "c"', shift right, curve={height=12pt}, from=1-1, to=1-3]
        \arrow[""{name=2, anchor=center, inner sep=0}, "b"{description}, from=1-3, to=1-1]
        \arrow[shift left=6, between={0.2}{0.8}, dotted, no head, from=2, to=0]
        \arrow[shift left=6, between={0.2}{0.8}, dotted, no head, from=2, to=1]
     \end{tikzcd}\]
  \end{minipage}
  \hfill
  \begin{minipage}{0.25\textwidth}
    \centering
    \resizebox{\linewidth}{!}%
    {%
      \begin{tikzpicture}[
          scale=1.5,
          ->-/.style={decoration={
              markings,
              mark=at position #1 with {\arrow[scale=1.3,>=stealth]{>}}
            },postaction={decorate}},
          -<-/.style={decoration={
              markings,
              mark=at position #1 with {\arrow[scale=1.3,>=stealth]{<}}
            },postaction={decorate}}
      ]

            \draw[thick] (0, 2) ellipse (1.6 and 0.3);
            
            \draw[thick] (-1.6, 2) .. controls (-1.8, 0.5) and (-1.2, -2) .. (0, -2.2) 
                          .. controls (1.2, -2) and (1.8, 0.5) .. (1.6, 2);
          
            \filldraw (1.25, 1.8) circle (1.5pt);
          
            \draw[thick] (0.2, 0.3) .. controls (0.45, -0.1) and (0.45, -0.6) .. (0.1, -1.0);
            \draw[thick] (0.25, 0.1) .. controls (0.05, -0.2) and (0.05, -0.5) .. (0.25, -0.8);

          \draw[very thick, cyan, dashed] (-1.4, -0.35) arc (180:0:0.74cm and 0.2cm);
          \draw[very thick, cyan, ->-=0.55] (0.1, -0.35) arc (360:180:0.75cm and 0.2cm);

          \draw[very thick, orange, -<-=0.42] (0,-0.2) ellipse (1.0cm and 1.4cm);

      \end{tikzpicture}  
    }
  \end{minipage}
  \hfill
  \begin{minipage}{0.24\textwidth}
    \centering
    \resizebox{\linewidth}{!}%
    {
      \begin{tikzpicture}[
        scale=1.5,
        ->-/.style={decoration={
            markings,
            mark=at position #1 with {\arrow[scale=1.3,>=stealth]{>}}
          },postaction={decorate}},
        -<-/.style={decoration={
            markings,
            mark=at position #1 with {\arrow[scale=1.3,>=stealth]{<}}
          },postaction={decorate}}
      ]

          \draw[line width=1.2pt, violet, dashed]
              (-1.4, -0.2) to[out=15, in=170] (0.1, -0.2);

          \draw[line width=1.2pt, violet, ->-=0.55]
              (0.1, -0.2) to[out=210, in=180] (0, 1.2)
              to[out=0, in=90] (1.2, 0)
              to[out=-90, in=0] (0, -1.2)
              to[out=180, in=0] (-1.4, -0.2);

              \draw[thick] (0, 2) ellipse (1.6 and 0.3);
              
              \draw[thick] (-1.6, 2) .. controls (-1.8, 0.5) and (-1.2, -2) .. (0, -2.2) 
                            .. controls (1.2, -2) and (1.8, 0.5) .. (1.6, 2);
            
              \filldraw (1.25, 1.8) circle (1.5pt);
            
              \draw[thick] (0.2, 0.3) .. controls (0.45, -0.1) and (0.45, -0.6) .. (0.1, -1.0);
              \draw[thick] (0.25, 0.1) .. controls (0.05, -0.2) and (0.05, -0.5) .. (0.25, -0.8);

      \end{tikzpicture}    
    }
  \end{minipage}
  \hfill

\caption{A generating family for the punctured torus with one stop, and the 
quiver associated to this generating familiy. Then two intersecting closed 
curves, and the resolution of their intersection point, which corresponds to 
the cone of that point.}
\label{fig:PuncturedTorusExample}
\end{figure}

The Fukaya category of the punctured torus is generated by two curves, 
$A$ and $B$, which cut the surface into a pentagon (whose fifth side is the 
boundary edge) as shown in Figure~\ref{fig:MarkedSurfaces}. As we see in 
Figure~\ref{fig:PuncturedTorusExample}, these curves have three interesting 
angles between them, giving us the quiver as shown. One can use 
the orientations shown to generate the $\mathbb{Z}/2$-graded Fukaya category 
or one can use a line field to build a $\mathbb{Z}$-graded Fukaya category 
where the degrees are also $|a|=|c|=1$ and $|b|=0$.

To build the two closed curves shown in the third subfigure of 
Figure~\ref{fig:PuncturedTorusExample},
we must glue $A$ to itself along the map $\lambda ba$ and $B$ to itself along the map 
$\mu cb$, respectively. Note that these curves are closed, so have a rank $1$ 
local system given by the invertible elements $\lambda$ and $\mu$. 
We see that there's a degree $1$ intersection point from the blue curve to 
the orange curve, and we can compute the algebraic representative 
of this by computing $\text{Hom}^\bullet((A,\lambda ba), (B, \mu cb))$. 
This is spanned by $\langle cba \rangle_0 \oplus \langle a, c \rangle_1$
(where the subscripts indicate the $\mathbb{Z}/2$-graded homological degree%
\footnote{If you used a line field grading to get the same quiver, now viewed 
as a $\mathbb{Z}$-graded algebra, then $cba$ would be placed in degree $2$ 
rather than degree $0$}). As usual the differential on $\text{Hom}^\bullet$ 
is given by $\delta f = \delta_B \circ f - (-1)^{|f|} f \circ \delta_A$, 
which for us means that $\delta(cba) = 0$, $\delta a = \mu cba$, and 
$\delta c = \lambda cba$. 
Then $\text{Ext}^\bullet = \text{ker}(\delta) \big / \text{im}(\delta)$
and in particular the grade $1$ piece $\text{Ext}^1$ is spanned by 
$a - \frac{\mu}{\lambda} c$, making this the algebraic representative of the 
unique intersection point. We can compute its cone by deforming the derivative 
on $(A \oplus B, \lambda ba^{11} + \mu cb^{22})$ to 
$\text{Cone}(a - \frac{\mu}{\lambda} c) = (A \oplus B, \lambda ba^{11} + 
\mu cb^{22} + (a- \frac{\mu}{\lambda} c)^{21})$. Then after a change of basis 
(read: a conjugation of this new $2 \times 2$ differential matrix) we get 
$(A \oplus B, a - \frac{\mu}{\lambda} c)$, which we recognize as the resolved 
curve shown in the last subfigure of Figure~\ref{fig:PuncturedTorusExample}.

The punctured torus is a great place to test conjectures because it's 
so small! In practice, computing $\text{Hom}^\bullet(\alpha,\beta)$ 
can be quite annoying since one has to look at all the homs from every 
indecomposable in $\alpha$ to every indecomposable in $\beta$ and then 
take cohomology. However, since in this example there are only two 
indecomposables $A$ and $B$ with small homsets, the computations can stay 
manageable by hand. The author has written some sagemath 
code\footnote{which she will hopefully make public soon} for computing with 
the punctured torus that even brings quite large exmaples into reach.

\subsection{(Generalized) Pairs of Pants}

\begin{figure}[htbp]
  \centering

  \begin{minipage}[c]{0.22\textwidth}
    \centering
    \resizebox{\linewidth}{!}%
    {%
      \begin{tikzpicture}[x=1cm, y=1cm, line width=1.2pt, >=Stealth]
          \draw (0,0) ellipse (3.8cm and 2.5cm);

          \filldraw[black] (0, 2.5) circle (0.12cm);

          \draw (-1.6, 0.5) circle (0.35cm);
          \draw[->] ({-1.6 + 0.6*cos(-80)}, {0.5 + 0.6*sin(-80)}) arc (-80:260:0.6cm);
          \node[above] at (-1.6, 1.1) {\Large $c$};

          \draw (1.6, 0.5) circle (0.35cm);
          \draw[->] ({1.6 + 0.6*cos(-80)}, {0.5 + 0.6*sin(-80)}) arc (-80:260:0.6cm);
          \node[above] at (1.6, 1.1) {\Large $a$};

          \draw[cyan] (-1.6, 0.15) -- (-1.6, -2.26);
          \draw[cyan] (1.6, 0.15) -- (1.6, -2.26);

          \node[text=cyan] at (-1.9, -2.7) {\Huge $A$};
          \node[text=cyan] at (1.9, -2.7) {\Huge $B$};

          \draw[->] (1.4, -1.8) to[bend left=4] (-1.4, -1.75);
          \node[above] at (0, -1.75) {\Large $b$};
      \end{tikzpicture}
    }
  \end{minipage}%
  \begin{minipage}[c]{0.22\textwidth}
    \centering
    $\begin{tikzcd}
      {} && {} \\
      A && B
      \arrow[""{name=0, anchor=center, inner sep=0}, draw=none, from=2-1, to=1-1]
      \arrow["c", from=2-1, to=2-1, loop, in=55, out=125, distance=10mm]
      \arrow[""{name=1, anchor=center, inner sep=0}, "b", from=2-1, to=2-3]
      \arrow[""{name=2, anchor=center, inner sep=0}, draw=none, from=2-3, to=1-3]
      \arrow["a", from=2-3, to=2-3, loop, in=55, out=125, distance=10mm]
      \arrow[shift left=3, curve={height=6pt}, between={0.2}{0.6}, dotted, no head, from=1, to=0]
      \arrow[shift right=3, curve={height=-6pt}, between={0.2}{0.6}, dotted, no head, from=1, to=2]
    \end{tikzcd}$
  \end{minipage}%
  \hfill
  \begin{minipage}[c]{0.24\textwidth}
    \centering
    \resizebox{\linewidth}{!}%
    {%
      \begin{tikzpicture}
          \draw[line width=1.3pt] (0,0) ellipse (4cm and 2.2cm);
          
          \draw[line width=1.3pt] (-1.6,0) circle (0.5cm);
          \draw[line width=1.3pt] (1.6,0) circle (0.5cm);
          
          \draw[color=cyan, line width=1.8pt] (-1.6,-0.5) -- (-1.6,-2.02);
          \node[color=cyan, below, font=\Large, yshift=-2mm] at (-1.6,-2.02) {$A$};
          
          \draw[color=cyan, line width=1.8pt] (1.6,-0.5) -- (1.6,-2.02);
          \node[color=cyan, font=\Large, yshift=-2mm] at (1.6,-2.02) [below] {$B$};
          
          \draw[color=cyan, line width=1.8pt] (-1.1,0) 
              to[out=20, in=180] (1.1, 0.75) 
              to[out=0, in=90] (2.35, 0.2) 
              to[out=270, in=0] (2.1, 0);
          \node[color=cyan, font=\Large] at (0.2, 0.2) {$C$};
          
          \draw[color=cyan, line width=1.8pt] (-2.1, 0) 
              to[out=150, in=180] (-1.0, 1.4) 
              -- (1.0, 1.4) 
              to[out=0, in=90] (3.2, 0) 
              to[out=270, in=330] (1.95, -0.35);
          \node[color=cyan, font=\Large] at (-0.2, 1.05) {$D$};
          
          \filldraw[fill=black, draw=black, line width=0.8pt] (0,2.2) circle (0.12cm);
          \filldraw[fill=black, draw=black, line width=0.8pt] (-1.6,0.5) circle (0.12cm);
          \filldraw[fill=black, draw=black, line width=0.8pt] (1.6,0.5) ellipse (0.12cm);

          \draw[-stealth, line width=1pt, <-] (-1.6, -0.7) to[out=180, in=270] (-2.3, 0.15);
          \node at (-2, -0.8) {\large $c_1$};
          
          \draw[-stealth, line width=1pt] (-1.6, -0.7) to[out=0, in=270] (-0.9, 0.1);
          \node at (-1.2, -0.8) {\large $c_2$};
          
          \draw[stealth-, line width=1pt] (-1.6, -1.6) -- (1.6, -1.6) node[midway, above] {\large $b$};
          
          \draw[-stealth, line width=1pt] (1.6, -0.8) to[out=0, in=270] (2.2, -0.45);
          \node at (2.15, -1) {\large $a_1$};
          
          \draw[-stealth, line width=1pt] (2.2,-0.45) to[out=70, in=270] (2.3, 0.1);
          \node at (2.7, 0) {\large $a_2$};
      \end{tikzpicture}
    }
  \end{minipage}%
  \begin{minipage}[c]{0.22\textwidth}
    \centering
    $\begin{tikzcd}
      C && D \\
      A && B
      \arrow[from=1-1, to=1-3]
      \arrow[""{name=0, anchor=center, inner sep=0}, from=1-1, to=2-1]
      \arrow[""{name=1, anchor=center, inner sep=0}, from=1-3, to=2-3]
      \arrow[""{name=2, anchor=center, inner sep=0}, from=2-1, to=1-3]
      \arrow[""{name=3, anchor=center, inner sep=0}, "b"', from=2-1, to=2-3]
      \arrow[shift left=2, between={0}{0.6}, dotted, no head, from=0, to=3]
      \arrow[shift left=3, between={0.4}{1}, dotted, no head, from=2, to=1]
    \end{tikzcd}$
  \end{minipage}

  \caption{Two different pairs of pants and their associated quivers. 
  With these, curves are allowed to end on the pant legs as well as the 
  waist. One can also consider a pair of pants where curves are not allowed 
  to end on the legs, as shown in Figure~\ref{fig:ThreeLegPairOfPants}. Note 
  that in the stopless version we have countably many self maps $c^n : A \to A$
  and $a^n : B \to B$ because we can wrap around the pant leg countably many 
  times, while in the case where every boundary has a stop our endomorphism 
  algebra is finite dimensional.}
  \label{fig:PantsVariations}
\end{figure}

\begin{figure}
  \centering
  \begin{minipage}{0.4\textwidth}
    \centering
    \resizebox{\linewidth}{!}%
    {%
      \begin{tikzpicture}[>=Stealth, scale=1.5]

        \draw[thick, draw=black!80] (0, 0) ellipse (3.4 and 1.3);

        \foreach \x in {-1.9, -0.5, 0.5, 1.9} {
          \node[text=black!90] at (\x, 0) {\Large $\times$};
        }

        \draw[cyan, line width=1.3pt, 
              decoration={markings, mark=at position 0.25 with {\arrow{>}}
              }, 
              postaction={decorate}] 
          (0,0) ellipse (1.5 and 0.75);

        \draw[orange, line width=1.3pt, 
              decoration={markings, 
                          mark=at position 0.25 with {\arrow{>}},
                          mark=at position 0.75 with {\arrow{>}}
                          }, 
              postaction={decorate}] 
          plot[domain=0:360, samples=120] ({2.4*sin(\x)}, {0.8*sin(2*\x)});

        \fill[black!80] (0,0) circle (2pt) node[below=5pt] {$r$};
        \fill[black!80] (-0.9487, 0.5809) circle (2pt) node[above=5pt] {$p$};
        \fill[black!80] (0.9487, 0.5809) circle (2pt) node[above=5pt] {$q$};

        \fill[black!80] (0, 1.3) circle (2pt);

      \end{tikzpicture}
    }
  \end{minipage}
  \hfill
  \begin{minipage}{0.25\textwidth}
    \centering
    \resizebox{\linewidth}{!}%
    {%
      \begin{tikzpicture}[
          ->-/.style={
              decoration={
                  markings,
                  mark=at position #1 with {\arrow[scale=1.3]{Stealth}}
              },
              postaction={decorate}
          },
          cross/.pic={
              \draw[thick, black] (-0.12, -0.12) -- (0.12, 0.12);
              \draw[thick, black] (-0.12, 0.12) -- (0.12, -0.12);
          }
      ]
          \draw[thick, black] (0, 0) ellipse (4.7cm and 2.2cm);
          \fill[black] (0, 2.2) circle (2pt);
          \pic at (-3.2, 0) {cross};
          \pic at (-1.2, 0) {cross};
          \pic at ( 1.2, 0) {cross};
          \pic at ( 3.2, 0) {cross};
          \draw[color=orange, thick, ->-=0.125]
              (0,0) .. controls (-0.9,  1.4) and (-2.5,  1.4) .. (-2.5, 0)
                    .. controls (-2.5, -1.4) and (-1,   -2.3) .. (2.5,  0)
                    .. controls ( 3.6,  0.8) and ( 4.3,  0.8) .. ( 4.3, 0)
                    .. controls ( 4.3, -1.4) and ( 0.9, -1.4) .. (0,    0);

          \draw[thick, cyan, ->-=0.4] (-2.2, 0) ellipse (1.6 and 0.7);
      \end{tikzpicture}
    }
  \end{minipage}
  \begin{minipage}{0.25\textwidth}
    \centering
    \resizebox{\linewidth}{!}%
    {%
      \begin{tikzpicture}[
          ->-/.style={
              decoration={
                  markings,
                  mark=at position #1 with {\arrow[scale=1.3]{Stealth}}
              },
              postaction={decorate}
          },
          cross/.pic={
              \draw[thick, black] (-0.12, -0.12) -- (0.12, 0.12);
              \draw[thick, black] (-0.12, 0.12) -- (0.12, -0.12);
          }
      ]
          \draw[thick, black] (0, 0) ellipse (4.7cm and 2.2cm);
          \fill[black] (0, 2.2) circle (2pt);
          \pic at (-3.2, 0) {cross};
          \pic at (-1.2, 0) {cross};
          \pic at ( 1.2, 0) {cross};
          \pic at ( 3.2, 0) {cross};
          \draw[color=violet, thick, ->-=0.125]
              (0,0) .. controls (-0.9,  1.4) and (-4.3,  1.4) .. (-4.3, 0)
                    .. controls (-4.3, -1.4) and (-1,   -2.3) .. (-0.8, 0)
                    .. controls ( -0.8,  0.8) and ( -2.3,  0.8) .. (-2.5, 0)
                    .. controls ( -2.7,  -1.8) and ( 2.2,  -0.7) .. ( 2.5, 0)
                    .. controls ( 3.4,  1.6) and ( 4.3,  0.8) .. ( 4.3, 0)
                    .. controls ( 4.3, -1.4) and ( 0.9, -1.4) .. (0,    0);

      \end{tikzpicture}
    }
  \end{minipage}
  \caption{The presence of the triangle $pqr$ ensures that $\text{Cone}(p+q)$
  will not be curve one gets by resolving first $p$ and then $q$
  (shown in the middle subfigure), but will instead be the resolution of 
  just $q$ (shown in the third subfigure).}
  \label{fig:BigPantsTriangle}
\end{figure}

Generalized pairs of pants, by which we mean surfaces of genus $0$, are 
very useful test beds for computation and conjecture. As special cases we 
recover disks with stops (and thus the $A_n$-quiver)
as in Figure~\ref{fig:AnQuiver}, and the annulus with stops (and thus the 
Kronecker quiver) as in Figure~\ref{fig:Monodromy}. In this setting one can 
clearly see the distinction between punctures, boundary components without 
stops, and boundary components with stops, as in 
Figure~\ref{fig:ThreeLegPairOfPants} and Figure~\ref{fig:PantsVariations}.

Moreover, since the genus $0$ case is easier to draw, it's easier to count 
holomorphic triangles that impact the algebra. Indeed, from our previous 
discussion one might expect that when curves intersect multiple times
(say, in degree $1$ points $p_1, \ldots, p_n$) that 
$\text{Cone}(\sum_{i=1}^n p_i)$ should be the resolution of each of the $n$ 
intersections points. This will be true provided the curves do not bound 
any polygons, but in general the situation is more complicated.

Consider, for example, the situation in Figure~\ref{fig:BigPantsTriangle}
where $p$ and $q$ are degree $1$ maps from the blue curve to the orange 
curve. Then in notation as in Figure~\ref{fig:ThreeLegPairOfPants} with a 
fourth curve $D$ and arrows $a$ through $g$, we see that the blue curve 
is given by the dg-module $ X = (C \oplus B, d^{21} + \lambda edc^{21})$ and the 
orange curve is 
$Y = (B \oplus A \oplus C[1] \oplus D[1], f^{21} + gfed^{23} + b^{34} + \mu dcba^{14})$.
Then one computes $\text{Hom}^\bullet(X,Y)$ by using the fact that 
$\text{Hom}^\bullet(\bigoplus P_i, \bigoplus Q_j) \cong \bigoplus_{i,j} 
\text{Hom}^\bullet(P_i, Q_j)$ and using the usual formula for the 
differential $\delta f = \delta_Y \circ f - (-1)^{|f|} f \circ \delta_X$. 
This gives us $\text{Hom}^\bullet(X,Y)$ as a $16$-dimensional vector space 
with $\delta$ a $16 \times 16$ matrix (which is thankfully quite sparse). 
One can than compute\footnote{using a computer algebra system, of course} 
the cohomology of $\delta$ to see that 
$\text{Ext}^1(X,Y)$ is spanned by two generators (each of which is a 
$4 \times 2$ matrix): $d^{11}$ (which is cohomologous to $f^{22} - \lambda edc^{11}$)
and $gfe^{22} - \text{Id}_C^{31}$, which represent the intersection points 
$p$ and $q$, respectively. It's a nice exercise to check this for yourself,
and to find algebraic representatives for the two degree $0$ intersection 
points as well.

Then one can compute $\text{Cone}(\alpha p + \beta q)$ to be the deformation of 
$X \oplus Y$ where we add $\alpha p + \beta q$ to $\delta_X + \delta_Y$, 
which gives us the dg-module 
$(C \oplus B) \oplus (B \oplus A \oplus C[1] \oplus D[1])$ with differential
$(d^{21} + \lambda edc^{21}) + (f^{43} + gfed^{45} + b^{56} + \mu dcba^{36})
+ (\alpha f^{42} - \alpha \lambda edc^{31}) + 
(\beta gfe^{42} - \beta \text{Id}_C^{51})$. 
After a suitable change of basis (again, conjugation by an appropriate 
matrix) this is the same dg-module as 
$(C \oplus C[1]) \oplus (B \oplus D[1]) \oplus (A \oplus B)$ with some 
differential. Then using that $(C \oplus C[1], \text{Id}_C^{21})$ is 
quasi-isomorphic to $0$ we see this is the same dg-module as 
$(B \oplus D[1]) \oplus (A \oplus B)$ with differential 
$(\mu dcba - \lambda \frac{\alpha}{\beta} edcb)^{12} +
(f + \frac{\alpha}{\beta} gfe)^{34} + 
(f^{31} + \lambda \frac{\alpha}{\beta} edcb^{42})$. If this differential 
only had the first two parenthesized terms then it would be the multicurve 
one gets by resolving both $p$ and $q$, as shown in the middle subfigure 
of Figure~\ref{fig:BigPantsTriangle}. However, the third parenthesized 
term is exactly the algebraic representative for the degree $1$ map 
from the blue curve to the orange curve in this middle subfigure! Resolving 
this extra crossing produces the purple curve in the third subfigure, 
which is thus the correct geometric representative of 
$\text{Cone}(\alpha p + \beta q)$. Note that this purple curve is also 
what one gets by simply resolving $q$, and this is not an accident! 
Indeed, if we resolve $q$ first then $p$ and $r$ bound a bigon which 
can be removed by an isotopy\footnote{Or said another way, this bigon witnesses
$p$ as the differential of $r$ in the Floer theory}. 
Thus after resolving $q$ the intersection point $p$ goes away! 
Algebraically this is beacuse $p$ becomes nullhomologous,
with $\delta r = p$ (perhaps up to an invertible factor), and the end result 
is that deforming by $q$ then by $p$ is the same thing as simply deforming 
by $q$.

\newpage
\printbibliography

\end{document}